\documentclass[12pt,letterpaper]{article}
\usepackage{graphicx}
\usepackage{amsmath}
\usepackage{amsfonts}
\usepackage{amsthm}
\usepackage{hyperref}

\begin{document}

\title{Tricky Arithmetic}

\author{Tanya Khovanova}

\maketitle

\abstract{This article is an expanded version of my talk at the Gathering for Gardner, 2012.}

%%%%%%%%%%%%%%%%%%%%%%%%%%%%%%%%%%%%%%%%%%%%%%%%%%%%%%%%%%%%%%%%%%%%%%%%%%%%%%%%
\section{Introduction}

\begin{quote}
Three horses are galloping at 27 miles per hour. What is the speed of one horse?
\end{quote}

If you answered 9 to the question above, you're wrong. Relax, stop rushing through this problem, and think again.

\section{Ode to trick problems}

I work with high school youth as a coach for math competitions. From time to time I give them trick problems. Many new students get caught, for example, answering 9 to the previous riddle. This may seem surprising, given that I work with the highest achievers in one of the best schools in Massachusetts.

I have a theory for why this happens. Few kids are really taught to think in school. Instead, they are taught to use templates. At a rushed first glance the problem above sounds like division, so they divide.

Once my students get caught for the first time, they start paying attention. After that, they think before answering. No one likes to be tricked, which is why trick problems are such a great motivator for students to pay attention and to really think.

\section{The collection of the problems}

As trick problems have great value, I collect them and put them on my website \cite{webpagetrick}. I also write about them for my blog \cite{blogtrick}. Here are a few gems from my collection:

\begin{enumerate}
\item A stick has two ends. If you cut off one end, how many ends will the stick have left? 
\item Anna had two sons. One son grew up and moved away. How many sons does Anna have now? 
\item A square has four corners. If we cut one corner off, how many corners will the remaining figure have? 
\item At a farmer's market you stop by an apple stand, where you see 20 beautiful apples. You buy 5. How many apples do you have? 
\item Mrs. Fullhouse has 2 sons, 3 daughters, 2 cats and 1 dog. How many children does she have? 
\item My dining room chandelier has 5 light bulbs in it. During a storm two of them went out. How many light bulbs are in the chandelier now? 
\item My dog Fudge likes books. In the morning he brought two books to his corner and three more books in the evening. How many books will he read tonight? 
\item There were five bowls full of candy on the table. Mike ate one bowl of candy and Sarah ate two. How many bowls are there on the table now?
\item Peter had ten cows. All but nine died. How many cows are left? 
\item A patient needs to get three shots. There is a break of 30 minutes between shots. Assuming the shots themselves are instantaneous, how much time will the procedure take? 
\item You are running a race and you pass the person who was running second. What place are you now? 
\item A caterpillar wants to see the world and decides to climb a 12-meter pole. It starts every morning and climbs 4 meters in half a day. Then it falls asleep for the second half of the day, during which time it slips 3 meters down. How much time will it take the caterpillar to reach the top? 
\item Humans have 10 fingers on their hands. How many fingers are there on 10 hands? 
\item How many people are there in two pairs of twins, twice? 
\item It takes 3 minutes to boil 3 eggs. How long will it take to boil 5 eggs? 
\item Two friends went for a walk and found \$20. How much money would they have found if they were joined by two more friends? 
\item One hundred percent of the fish in a pond are goldfish. I take 10\% of the goldfish out of the pond. What percentage of goldfish does the pond now contain? 
\item The Davidsons have five sons. Each son has one sister. How many children are there in the family? 
\end{enumerate}

\section{The wrong answers}

Here are the wrong answers that are produced by people who use a template without thinking:

\begin{enumerate}
\item 1 end
\item 1 son
\item 3 corners
\item 15 apples
\item 8 children
\item 3 light-bulbs
\item 5 books
\item 2 bowls
\item 1 cow
\item One hour and a half
\item First
\item 12 days
\item 100 fingers
\item 16 people
\item 5 minutes
\item \$40
\item 90\% goldfish
\item 10 children
\end{enumerate}

\section{Rabbits}

I know that my students are doing something wrong when they start using a pencil or a calculator for this problem:

\begin{quote}
On average, rabbits start breeding when they are 3 months old and produce 4 offspring every month. If I put a day old rabbit in a cage for a year, how many offspring will it produce? 
\end{quote}

\section{Acknowledgments}

Thanks to all those who submitted trick problems to me. I am grateful to Sue Katz for helping me to edit this paper.

\end{document}